\documentclass{article}

\usepackage{amssymb,latexsym,amsmath}

\usepackage[pdftex]{graphicx}

\usepackage{hyperref}

\begin{document}

\hoffset-0.64cm
\voffset-2.14cm

\textheight22.8cm

\textwidth15.5cm

\pagestyle{plain}

\newcommand{\bfi}{\bfseries\itshape}

\makeatletter

\@addtoreset{figure}{section}

\def\thefigure{\thesection.\@arabic\c@figure}

\def\fps@figure{h, t}

\@addtoreset{table}{bsection}

\def\thetable{\thesection.\@arabic\c@table}

\def\fps@table{h, t}

\@addtoreset{equation}{section}

\def\theequation{\thesubsection.\arabic{equation}}

\makeatother

\newtheorem{theorem}{Theorem}[section]

\newtheorem{proposition}[theorem]{Proposition}

\newtheorem{lema}[theorem]{Lemma}

\newtheorem{corollary}[theorem]{Corollary}

\newtheorem{definition}[theorem]{Definition}

\newtheorem{remark}[theorem]{Remark}

\newtheorem{exempl}{Example}[section]

\newenvironment{exemplu}{\begin{exempl}  \em}{\hfill $\square$

\end{exempl}}

\newcommand{\comment}[1]{\par\noindent{\raggedright\texttt{#1}

\par\marginpar{\textsc{Comment}}}}

\newcommand{\todo}[1]{\vspace{5 mm}\par \noindent \marginpar{\textsc{ToDo}}\framebox{\begin{minipage}[c]{0.95 \textwidth}

\tt #1 \end{minipage}}\vspace{5 mm}\par}

\newcommand{\ea}{\mbox{{\bf a}}}

\newcommand{\eu}{\mbox{{\bf u}}}

\newcommand{\ueu}{\underline{\eu}}

\newcommand{\ueo}{\overline{u}}

\newcommand{\oeu}{\overline{\eu}}

\newcommand{\ew}{\mbox{{\bf w}}}

\newcommand{\ef}{\mbox{{\bf f}}}

\newcommand{\eF}{\mbox{{\bf F}}}

\newcommand{\eC}{\mbox{{\bf C}}}

\newcommand{\en}{\mbox{{\bf n}}}

\newcommand{\eT}{\mbox{{\bf T}}}

\newcommand{\eL}{\mbox{{\bf L}}}

\newcommand{\eR}{\mbox{{\bf R}}}

\newcommand{\eV}{\mbox{{\bf V}}}

\newcommand{\eU}{\mbox{{\bf U}}}

\newcommand{\ev}{\mbox{{\bf v}}}

\newcommand{\eve}{\mbox{{\bf e}}}

\newcommand{\uev}{\underline{\ev}}

\newcommand{\eY}{\mbox{{\bf Y}}}

\newcommand{\eK}{\mbox{{\bf K}}}

\newcommand{\eP}{\mbox{{\bf P}}}

\newcommand{\eS}{\mbox{{\bf S}}}

\newcommand{\eJ}{\mbox{{\bf J}}}

\newcommand{\eB}{\mbox{{\bf B}}}

\newcommand{\eH}{\mbox{{\bf H}}}

\newcommand{\leb}{\mathcal{ L}^{n}}

\newcommand{\eI}{\mathcal{ I}}

\newcommand{\eE}{\mathcal{ E}}

\newcommand{\hen}{\mathcal{H}^{n-1}}

\newcommand{\eBV}{\mbox{{\bf BV}}}

\newcommand{\eA}{\mbox{{\bf A}}}

\newcommand{\eSBV}{\mbox{{\bf SBV}}}

\newcommand{\eBD}{\mbox{{\bf BD}}}

\newcommand{\eSBD}{\mbox{{\bf SBD}}}

\newcommand{\ecs}{\mbox{{\bf X}}}

\newcommand{\eg}{\mbox{{\bf g}}}

\newcommand{\paromega}{\partial \Omega}

\newcommand{\gau}{\Gamma_{u}}

\newcommand{\gaf}{\Gamma_{f}}

\newcommand{\sig}{{\bf \sigma}}

\newcommand{\gac}{\Gamma_{\mbox{{\bf c}}}}

\newcommand{\deu}{\dot{\eu}}

\newcommand{\dueu}{\underline{\deu}}

\newcommand{\dev}{\dot{\ev}}

\newcommand{\duev}{\underline{\dev}}

\newcommand{\weak}{\stackrel{w}{\approx}}

\newcommand{\mild}{\stackrel{m}{\approx}}

\newcommand{\lrightarrow}{\stackrel{L}{\rightarrow}}

\newcommand{\rrightarrow}{\stackrel{R}{\rightarrow}}

\newcommand{\strong}{\stackrel{s}{\approx}}

\newcommand{\weakdown}{\rightharpoondown}

\newcommand{\opg}{\stackrel{\mathfrak{g}}{\cdot}}

\newcommand{\opunu}{\stackrel{1}{\cdot}}
\newcommand{\opdoi}{\stackrel{2}{\cdot}}

\newcommand{\opn}{\stackrel{\mathfrak{n}}{\cdot}}
\newcommand{\opx}{\stackrel{x}{\cdot}}

\newcommand{\tr}{\ \mbox{tr}}

\newcommand{\Ad}{\ \mbox{Ad}}

\newcommand{\ad}{\ \mbox{ad}}

\renewcommand{\contentsname}{ }

\title{Maps of metric spaces}

\author{Marius Buliga \\ 
\\
Institute of Mathematics, Romanian Academy \\
P.O. BOX 1-764, RO 014700\\
Bucure\c sti, Romania\\
{\footnotesize Marius.Buliga@imar.ro}}

\date{This version:  09.06.2011}

\maketitle

\begin{abstract}
This is a pedagogical introduction covering maps of metric spaces, Gromov-Hausdorff
distance and its "physical" meaning, and dilation structures as a convenient
simplification of an exhaustive database of maps of a metric space into another. 
\end{abstract}










\section{Exploring space}
Suppose we send an explorer to make maps of  an unknown territory  
$X$. The explorer  wants to  record her discoveries on maps, or charts,  
done in the metric space $(Y,D)$.

\vspace{.5cm}

\centerline{\includegraphics[angle=270, width=0.5\textwidth]{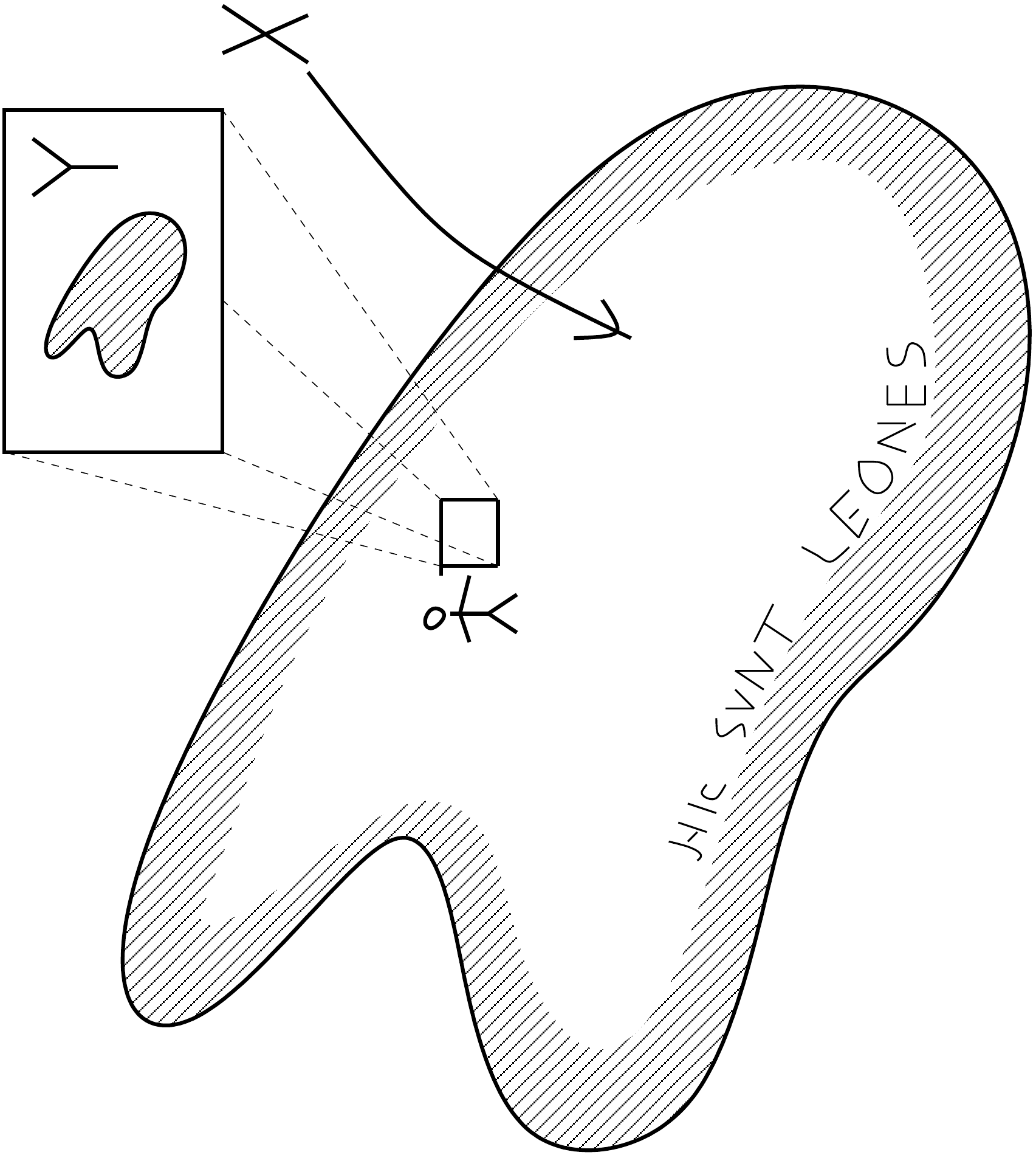}}

\vspace{.5cm}

I shall suppose that we can put a distance on the set $X$, that is a function 
$d: X \times X \rightarrow [0,+\infty)$ which satisfies the following 
requirement:  for any three points $x,y,z \in X$ there is a bijective 
correspondence with a triple A,B,C in the plane  such that the sizes (lengths) 
of AB, BC, AC are equal respectively with $d(x,y)$, $d(y,z)$, $d(z,x)$. 
Basically, we accept that we can represent in the plane any three points from the 
space $X$. An interpretation of the distance $d(x,y)$ is the following:
the explorer has a ruler and $d(x,y)$ is the numerical value shown by the ruler 
when streched between points marked with "$x$" and "$y$". (Then the explorer has
to use somehow these numbers in order to make a chart of $X$.)

\paragraph{How many maps are needed?}
 "Understanding"  the space $X$ (with respect 
to the choice of the "gauge" function $d$) into the terms of the more familiar 
space $Y$   means  making a chart $f:X \rightarrow Y$ 
of $X$ into $Y$ which is not deforming distances too much. Ideally, a 
perfect chart has to be Lipschitz, that is the distances between points in 
$X$ are transformed by the chart into distances between points in $Y$, with a 
precision independent of the scale: the chart  
$f$ is (bi)Lipschitz if there are positive numbers 
$c < C$ such that for any two points 
$x,y \in X$  
$$ c \, d(x,y) \, \leq \, D(f(x),f(y)) \, \leq \, C \, d(x,y)$$
This would be a very  good chart, because it would tell how $X$ is at all scales. 
There are two difficulties related to this model. First, it is impossible to 
make such a chart in practice. What we can do 
instead, is to sample the space $X$ (take a $\varepsilon$-dense subset of $X$ 
with respect to the distance $d$) and try to represent as good as possible this 
subspace in $Y$. Mathematically this is like asking for the chart function $f$ to have the following property: there are supplementary positive constants $a, A$ such that 
for any two points 
$x,y \in X$  
$$ c \, d(x,y) - a  \, \leq \, D(f(x),f(y)) \, \leq \, C \, d(x,y) + A$$
The second difficulty is that such a chart might not exist at all, from
mathematical reasons (there is no quasi-isometry between the metric spaces 
$(X,d)$ and $(Y,D)$). Such a chart exists of course if we want to make charts of 
regions with bounded distance, but remark that all details are erased at small
scale. The remedy would be to make better and better charts, at smaller and
smaller scales, eventually obtaining something resembling a road atlas, with 
charts of countries, regions, counties, cities, charts which have to be compatible 
one with another in a clear sense.

\section{From maps to dilation structures}

Imagine that  the metric space $(X,d)$ represents a territory. We want to make 
maps  of $(X,d)$ in the metric space $(Y,D)$ (a piece of paper, or a
scaled model). 

In fact, in order to understand the territory $(X,d)$, we need many maps, at
many scales. For any point $x \in X$ and any scale $\varepsilon > 0$ we shall 
make a map of a neighbourhood of $x$, ideally. In practice, a good knowledge 
of a territory amounts to have, for each of several scales $\displaystyle
\varepsilon_{1} > \varepsilon_{2} > ... > \varepsilon_{n}$ an atlas of maps 
of overlapping parts of  $X$ (which together form a cover of the territory $X$). 
All the maps from all the atlasses have to be compatible one with another. 

The ideal model of such a body of knowledge is embodied into the notion of 
a manifold. To have $X$ as a manifold over the model space $Y$ means exactly 
this. 

Examples from metric geometry (like sub-riemannian spaces) show that the 
manifold idea could be too rigid in some situations. We shall replace it 
with the idea of a dilation structure. 

We shall see that a dilation structure (the right generalization of a 
smooth space, like a manifold), represents an idealization 
of the more realistic situation of having at our disposal many maps, at many
scales, of the territory, with the property that the accuracy, precision and 
resolution of such maps, and of relative maps deduced from them, are controlled 
by the scale (as the scale goes to zero, to infinitesimal). 

There are two facts which I need to stress. First is that such a generalization
is necessary. Indeed, by looking at the large gallery of metric spaces which we
now know, the metric spaces with a manifold structure form a tiny and very very
particular class. Second is that we tend to take for granted the body of
knowledge represented by a manifold structure (or by a dilation structure). 
Think as an example at the manifold structure of the Earth. It is an
idealization of the collection of all cartographic maps of parts of the Earth. 
This is a huge data basis and it required a huge amount of time and energy in
order to be constructed. To know, understand the territory is a huge task, 
largely neglected. We "have" a manifold, "let $X$ be a manifold". And even if 
we do not doubt that the physical space (whatever that means) is a boring 
$\displaystyle \mathbb{R}^{3}$, it is nevertheless another task to determine 
with the best accuracy possible a certain point in that physical space, based on
the knowledge of the coordinates. For example GPS costs money and time to build 
and use. Or, it is rather easy to collide protons, but to understand and keep 
the territory fixed (more or less) with respect to the map, that is where most
of the effort goes. 

A model of such a map of $(X,d)$ in $(Y,D)$ is  a relation  
$\rho \subset X \times Y$,  a subset of 
a cartesian product $X \times Y$ of two sets. A particular type of relation 
is the graph of a function $\displaystyle f: X \rightarrow Y$, defined as the
relation 
$$\rho \, = \, \left\{ (x, f(x)) \mbox{ : } x \in X \right\}$$ 
but there are many relations which cannot be described as graphs of functions.


\centerline{\includegraphics[angle=270, width=0.5\textwidth]{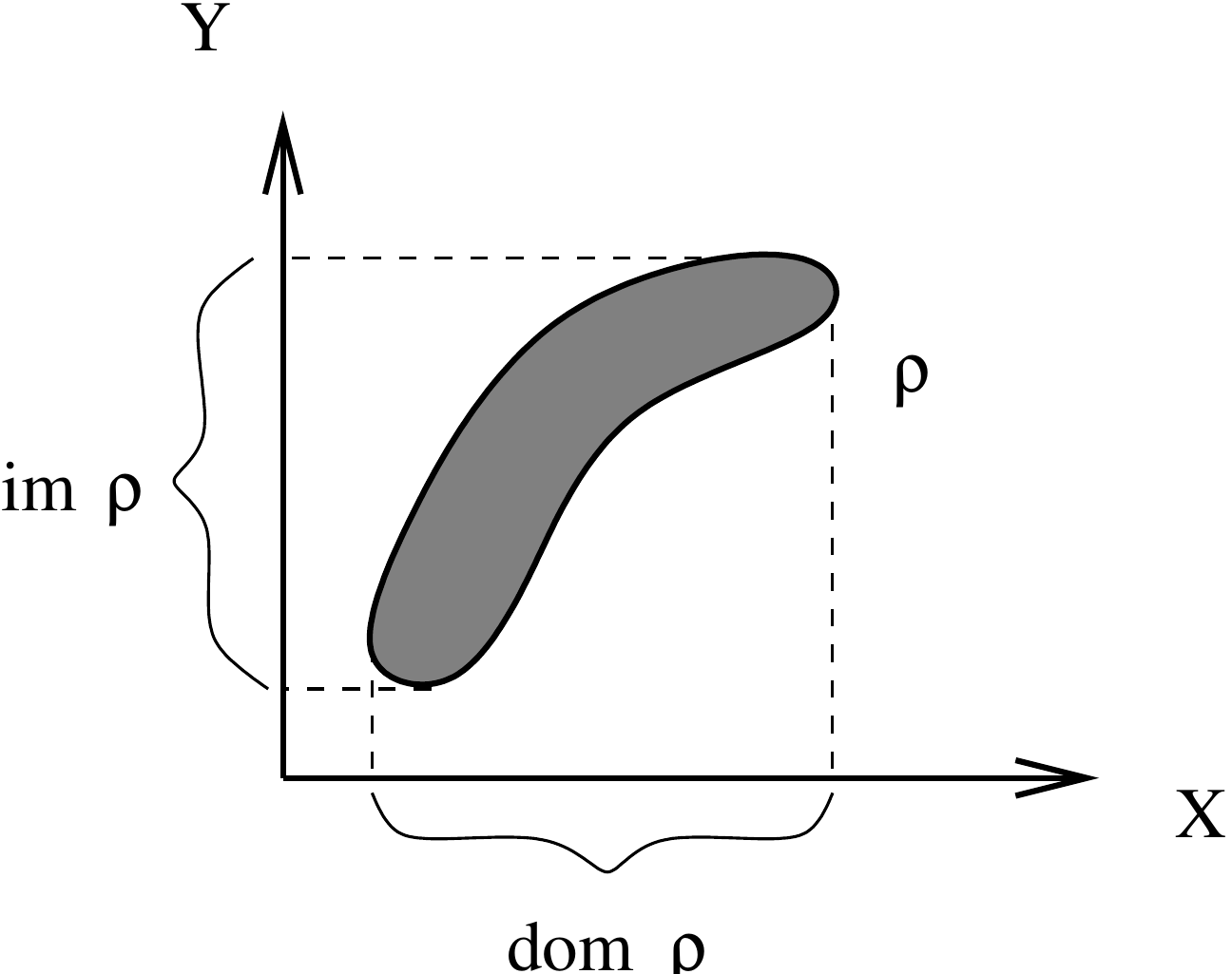}}


Imagine that pairs $(u,u') \in \rho$ are pairs 

\vspace{.5cm}

\centerline{(point in the space $X$, pixel in the "map space" Y)} 

\vspace{.5cm}

\noindent with the meaning that 
the point $u$ in $X$ is represented as the pixel $u'$ in $Y$. 

I don't suppose that there is a one-to-one correspondence between points in 
$X$ and pixels in $Y$, for various reasons, for example: due to repeated measurements 
there is no unique way to associate pixel to a point, or a point to a pixel. 
The relation $\rho$ represents  the cloud of pairs point-pixel which are 
compatible with all measurements. 

I shall use this model of a map for simplicity reasons. A better,  
more realistic model could be one using probability measures, but this model is 
sufficient for the needs of this paper. 

For a given map $\rho$ the point $x \in X$ in the space $X$  is associated
 the set of points $\left\{ y \in Y \mbox{ : } (x,y) \in \rho \right\}$ in the 
 "map space" $Y$. Similarly, to the "pixel" $y \in Y$ in the "map space" $Y$
  is associated the set of points 
  $\left\{ x \in X \mbox{ : } (x,y) \in \rho \right\}$ in the space $X$. 

\vspace{.5cm}

\centerline{\includegraphics[angle=270, width=0.5\textwidth]{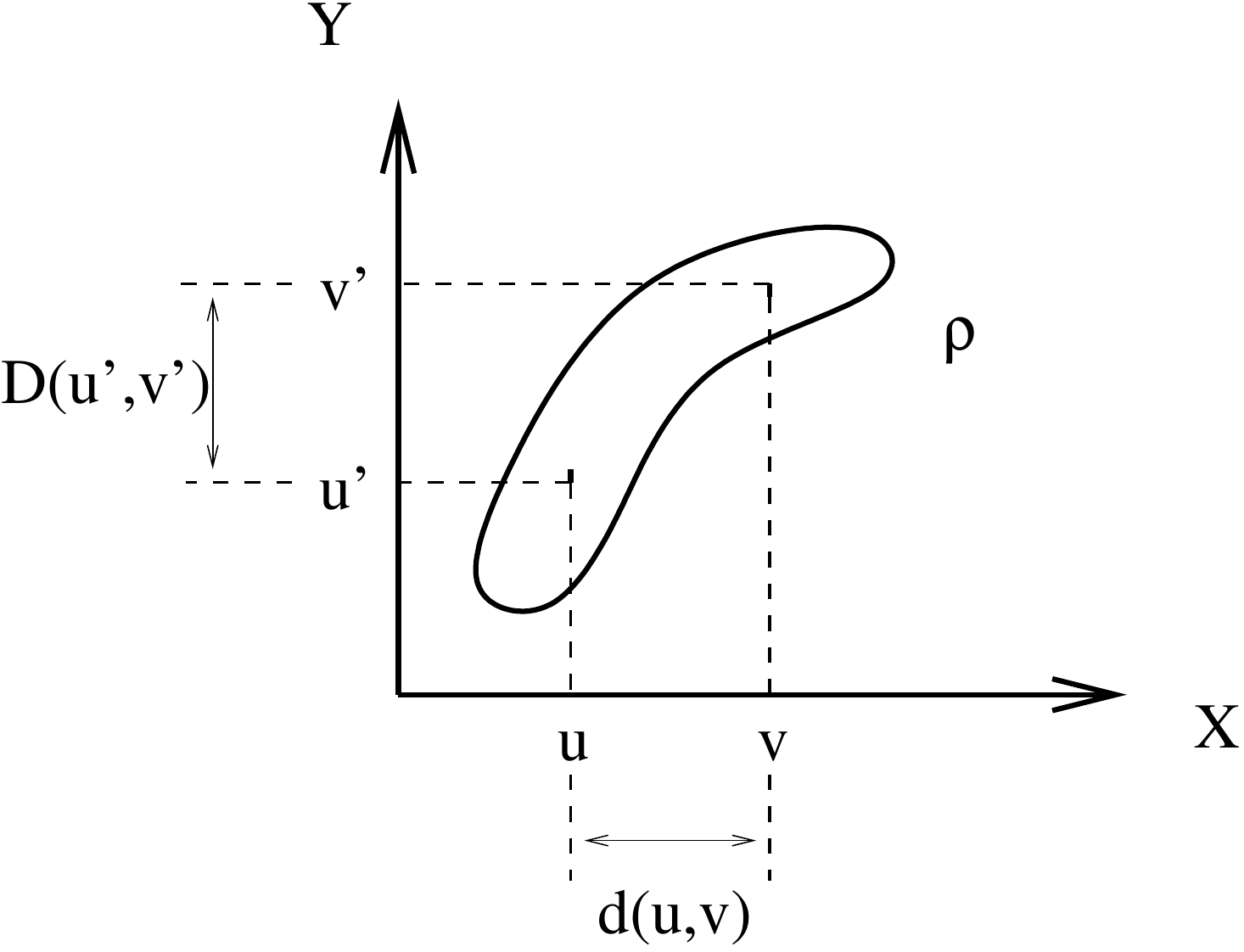}}

\vspace{.5cm}

A good map is one which does not distort distances too much. Specifically,
considering any two points $u,v \in X$ and any two pixels $u',v' \in Y$ which 
represent these points, i.e. $(u,u'), (v,v') \in \rho$, the distortion of 
distances between these points is measured by the number

$$\mid d(u,v) - D(u',v') \mid$$

\section{Accuracy, precision, resolution, Gromov-Hausdorff distance}

\paragraph{Notations concerning relations.}  
Even if relations are more general than (graphs of) functions, there is no harm to
use, if needed, a functional notation. For any relation  
$\rho  \subset X \times Y$  we shall write $\rho(x) = y$ or $\displaystyle 
\rho^{-1}(y) = x$ if $(x,y) \in \rho$. Therefore we may have $\rho(x) = y$ and 
$\rho(x) = y'$ with $y \not = y'$, if $(x,y) \in f$ and $(x,y') \in f$. In some
drawings, relations will be figured by a large arrow, as shown further.

\vspace{.5cm}

\centerline{\includegraphics[angle=270, width=0.5\textwidth]{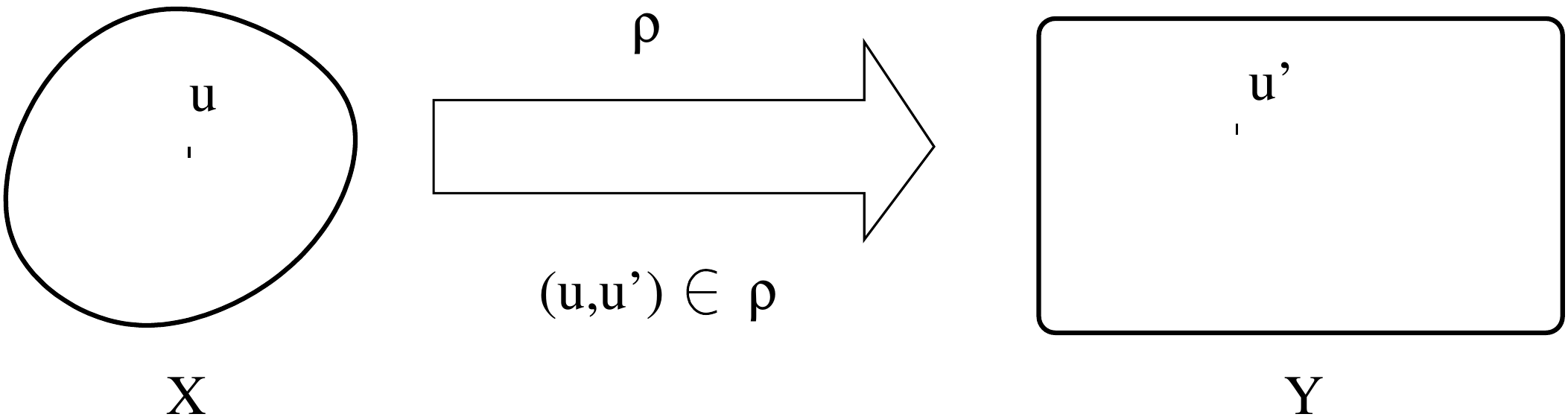}}

\vspace{.5cm}

The domain of the relation $\rho$ is the set $dom \ \rho \,  \subset X$ such that for any $x \in \, dom \ \rho$ 
there is $y \in Y$ with $\rho(x) = y$. The image of $\rho$ is the set of $im \ \rho \  \subset Y$ such that for any $y \in \, im \ \rho$ there is $x \in X$ with $\rho(x) =  y$.  
By convention, when we write that a statement $R(f(x), f(y), ...)$ is true, we mean that $R(x',y', ...)$ is 
true for any choice of $x', y', ...$, such that $(x,x'), (y,y'), ... \in f$. 

The inverse of a relation $\rho \subset X \times Y$ is the relation 
$$\rho^{-1} \subset Y \times X \, , \quad \rho^{-1} \, = \, \left\{ (u',u) \mbox{
: } (u,u') \in \rho \right\}$$ 
and if $\rho' \subset X \times Y$, $\rho" \subset Y \times Z$ are two relations,
their composition is defined as 
$$\rho = \rho" \circ \rho' \subset X \times Z $$ 
$$\rho \, = \, \left\{ (u,u") \in X \times Z \mbox{ : } \exists u'\in Y \, (u,u')
\in \rho' \, (u',u") \in \rho" \right\}$$

I shall use the following convenient notation: by $\mathcal{O}(\varepsilon)$ we mean a positive function such that $\displaystyle \lim_{\varepsilon \rightarrow 0} \mathcal{O}(\varepsilon) = 0$.

\vspace{.5cm}


In metrology, by definition, accuracy is \cite{metrology} 2.13 (3.5) "closeness of 
agreement between a measured quantity value and a true quantity value of a
measurand". (Measurement) precision is \cite{metrology} 2.15 "closeness of agreement between 
indications or measured quantity values obtained by replicate
measurements on the same or similar objects under specified conditions".
Resolution  is \cite{metrology} 2.15 "smallest change in a quantity being measured that
causes a perceptible change in the corresponding
indication".

For our model of a map, if $(u,u') \in \rho$ then $u'$ represent the measurement
of $u$. Moreover, because we see a map as a relation, the definition of 
the resolution can be restated as the supremum of distances between points 
in $X$ which are represented by the same pixel. Indeed, if the distance between 
two points in $X$ is bigger than this supremum then they cannot be represented 
by the same pixel.

\begin{definition}
Let $\rho \subset X \times Y$ be a relation which represents a map of 
$dom \ \rho \, \subset (X,d)$ into $ im \ \rho \, \subset (Y,D)$. To this map 
are associated three quantities: accuracy, precision and resolution. 

 The accuracy of $\rho$ is defined by: 
\begin{equation}
acc(\rho) \, = \, \sup \left\{ \mid D(y_{1}, y_{2}) - d(x_{1},x_{2}) \mid \mbox{
: } (x_{1},y_{1}) \in \rho \, , \, (x_{2},y_{2}) \in \rho \right\}
\label{acc1}
\end{equation}
 The resolution of $\rho$ at  $y \in im \ \rho$ is 
\begin{equation}
res(\rho)(y) \, = \, \sup \left\{ d(x_{1},x_{2}) \mbox{
: } (x_{1},y) \in \rho \, , \, (x_{2},y) \in \rho \right\}
\label{resy1}
\end{equation}
and the resolution of $\rho$ is given by: 
\begin{equation}
res(\rho) \, = \, \sup \left\{ res(\rho)(y) \mbox{
: } y \in \, im \ \rho \right\}
\label{res1}
\end{equation}
 The precision of $\rho$ at $x \in dom \ \rho$ is 
\begin{equation}
prec(\rho)(x) \, = \, \sup \left\{ D(y_{1},y_{2}) \mbox{
: } (x,y_{1}) \in \rho \, , \, (x,y_{2}) \in \rho \right\}
\label{precx1}
\end{equation}
and the precision of $\rho$ is given by: 
\begin{equation}
prec(\rho) \, = \, \sup \left\{ prec(\rho)(x) \mbox{
: } x \in \, dom \ \rho \right\}
\label{prec1}
\end{equation}
\label{defacc}
\end{definition}

After measuring (or using other means to deduce) the distances $d(x', x")$ between 
all pairs of points in $X$ (we may have several values for the distance $d(x',x")$), 
we try to represent the collection 
of these distances  in $(Y,D)$. 
When we make a map $\rho$ we are not really measuring the distances between 
all points in $X$, then representing them as accurately as possible in $Y$.

What we do is that  we consider a relation $\rho$, with domain $M = \,
dom(\rho)$ 
 which is $\varepsilon$-dense in $(X,d)$, then we perform a " cartographic generalization"\footnote{\url{http://en.wikipedia.org/wiki/Cartographic_generalization},
"Cartographic generalization is the method whereby information is selected and 
represented on a map in a way that adapts to the scale of the display medium 
of the map, not necessarily preserving all intricate geographical or other 
cartographic details.} of the relation $\rho$ to a relation $\displaystyle \bar{\rho}$, a
 map of $(X,d)$ in $(Y,D)$,  for example as in the following definition. 

\begin{definition}
A subset $M \subset X$ of a metric space $(X,d)$ is $\varepsilon$-dense in 
$X$ if for any $u \in X$ there is $x \in M$ such that $d(x,u) \leq \varepsilon$. 

Let $\rho \subset X \times Y$ be a relation such that $dom \ \rho$ is 
$\varepsilon$-dense in $(X,d)$ and $im \ \rho$  is $\mu$-dense in 
 $(Y,D)$. We define then $\displaystyle \bar{\rho} \subset X \times Y$ by:
 $\displaystyle (x,y) \in \bar{\rho}$ if there is $\displaystyle (x',y') \in
 \rho$ such that $\displaystyle d(x,x')\leq \varepsilon$ and $\displaystyle 
 D(y,y') \leq \mu$.
 \label{defgencart1}
\end{definition}

If $\rho$ is a relation as described in definition \ref{defgencart1} then 
accuracy $acc(\rho)$, $\varepsilon$ and $\mu$ control the precision 
$prec(\rho)$ and resolution $res(\rho)$. Moreover, the accuracy, precision and 
resolution of $\displaystyle \bar{\rho}$ are controlled by those of 
$\rho$ and $\varepsilon$, $\mu$, as well. This is explained in the next
proposition.

\begin{proposition}
Let $\rho$ and $\displaystyle \bar{\rho}$ be as described in definition
\ref{defgencart1}. Then: 
\begin{enumerate}
\item[(a)] $\displaystyle res(\rho) \, \leq \, acc(\rho)$, 
\item[(b)] $\displaystyle prec(\rho) \, \leq \, acc(\rho)$, 
\item[(c)] $\displaystyle res(\rho) + 2 \varepsilon \leq \, res(\bar{\rho}) \leq
\, acc(\rho) + 2(\varepsilon + \mu)$,
\item[(d)] $\displaystyle prec(\rho) + 2 \mu \leq \, prec(\bar{\rho}) \leq
\, acc(\rho) + 2(\varepsilon + \mu)$,
\item[(e)] $\displaystyle \mid acc(\bar{\rho}) - \, acc(\rho) \mid \leq
2(\varepsilon + \mu)$. 
\end{enumerate}
\label{propacc1}
\end{proposition}

\paragraph{Proof.} Remark that (a), (b) are immediate consequences of definition
\ref{defacc} and that 
(c) and (d) must have identical proofs, just by switching $\varepsilon$ with 
$\mu$ and $X$ with $Y$ respectively. I shall prove therefore (c) and (e). 

For proving (c), consider $y \in Y$. By definition of 
$\displaystyle \bar{\rho}$ we write 
$$\left\{ x \in X \mbox{ : } (x,y) \in \bar{\rho} \right\} \, = \, \bigcup_{(x',y')
\in \rho , y' \in \bar{B}(y,\mu)} \bar{B}(x',\varepsilon)$$
Therefore we get 
$$res(\bar{\rho})(y) \, \geq \, 2 \varepsilon + \sup \left\{ res(\rho)(y') 
\mbox{ : } y' \in \,
im(\rho) \cap \bar{B}(y,\mu) \right\} $$
By taking the supremum over all $y \in Y$ we obtain the inequality 
$$res(\rho) + 2 \varepsilon \leq \, res(\bar{\rho})$$
For the other inequality, let us consider $\displaystyle (x_{1},y), (x_{2},y)
\in \bar{\rho}$ and $\displaystyle (x_{1}', y_{1}'), (x_{2}', y_{2}') \in \rho$
such that $\displaystyle d(x_{1},x_{1}') \leq \varepsilon, d(x_{2},x_{2}') 
\leq \varepsilon, D(y_{1}',y) \leq \mu,  D(y_{2}',y) \leq \mu$. Then: 
$$d(x_{1},x_{2}) \leq 2 \varepsilon + d(x_{1}',x_{2}') \leq 2 \varepsilon + \, 
acc(\rho) + d(y_{1}',y_{2}') \leq 2 (\varepsilon + \mu) + \, 
acc(\rho)$$
Take now a supremum and arrive to the desired inequality. 

For the proof of (e)  let us consider for $i = 1,2$   
$\displaystyle (x_{i},y_{i}) \in \bar{\rho}, (x_{i}',y_{i}') \in \rho$ such 
that $\displaystyle d(x_{i}, x_{i}') \leq \varepsilon, D(y_{i},y_{i}') \leq
\mu$.   It is then  enough to take absolute values and transform 
the following equality  
$$d(x_{1},x_{2}) - D(y_{1},y_{2}) = d(x_{1},x_{2}) - d(x_{1}',x_{2}') + 
d(x_{1}',x_{2}') - D(y_{1}',y_{2}') + $$ 
$$+ D(y_{1}',y_{2}') - D(y_{1},y_{2})$$ 
into well chosen, but straightforward, inequalities. \hfill $\square$

\vspace{.5cm}

The  following definition of the  Gromov-Hausdorff distance for metric spaces is natural, 
owing to the fact that the accuracy (as defined in definition \ref{defacc}) controls the 
precision and resolution. 

\begin{definition}
Let $\displaystyle (X, d)$, $(Y,D)$,  be a pair of metric spaces and $\mu > 0$. 
We shall say that $\mu$ is admissible if  there is a relation 
$\displaystyle \rho \subset X \times Y$ such that  
$\displaystyle dom \ \rho = X$, $\displaystyle im \ \rho = Y$, and $acc(\rho) \leq \mu$.
The Gromov-Hausdorff distance  between $\displaystyle (X,d)$ and $\displaystyle  
(Y,D)$  is   the infimum of admissible numbers $\mu$. 
\label{defgh}
\end{definition}

As introduced in definition \ref{defgh}, the Gromov-Hausdorff  (GH) distance is not a true 
distance, because the GH distance between two isometric  metric spaces   
is equal to zero. In fact the GH distance induces a distance on isometry classes of 
compact metric spaces. 

The GH distance thus represents a lower bound on the accuracy of making maps of 
$(X,d)$ into $(Y,D)$. Surprising as it might seem, there are many examples of pairs of 
metric spaces with the property that the 
GH distance between any pair of closed balls from these spaces, 
considered with the distances properly rescaled, is greater than a strictly
positive number, independent of the choice of the balls. Simply put: {\it there 
are pairs of spaces $X$, $Y$ such that is impossible to make maps of parts 
of $X$ into $Y$ with arbitrarily small accuracy.} 

Any measurement is equivalent with making a map, say of $X$ (the territory 
of the phenomenon) into $Y$ (the map space of the laboratory). The possibility
that there might a physical difference (manifested as a strictly positive 
GH distance) between these two spaces, even if they both might be topologically 
the same (and with trivial topology, say of a $\displaystyle \mathbb{R}^{n}$),
 is ignored in physics, apparently. On one
side, there is no experimental way to confirm that a territory is the same 
at any scale (see the section dedicated to the notion of scale), but much of physical
explanations are based on differential calculus, which has as the most basic
assumption that locally and infinitesimally the territory is the same. 
On the other side the imposibility of making maps of the phase space of a 
quantum object into the macroscopic map space of the laboratory might be a 
manifestation of the fact that there is a difference (positive GH distance 
between maps of the territory realised with the help of physical phenomena) 
between "small" and "macroscopic" scale.

\section{Scale} 

 Let $\varepsilon >0$. A map of $(X,d)$ into $(Y,D)$, at scale 
$\varepsilon$ is a map of $\displaystyle (X, \frac{1}{\varepsilon}d)$ into $(Y,D)$. 
Indeed, if this map would have accuracy equal to $0$ then a value of a 
distance between points in $X$ equal to  $\displaystyle L$ would correspond to a value 
of the distance between the corresponding points on the map in $(Y,D)$ equal to
$\varepsilon \, L$. 

In cartography, maps  of the same territory done at smaller and smaller scales (smaller 
and smaller $\varepsilon$) must have the property that, at the same resolution, 
the  accuracy and precision (as defined in definition \ref{defacc}) have to become 
smaller and smaller. 

In mathematics, this could serve as the definition of the metric tangent space to 
a point in $(X,d)$, as seen in $(Y,D)$. 

\begin{definition}
We say that $(Y,D, y)$ ($y \in Y$) represents the (pointed unit ball in the) metric 
tangent space at $x \in X$ of $(X,d)$ if there exist a pair formed by: 
\begin{enumerate}
\item[-] a "zoom sequence", that is a sequence 
$$\displaystyle \varepsilon \in (0,1] \, \mapsto \,  \rho_{\varepsilon}^{x} \subset
(\bar{B}(x,\varepsilon), \frac{1}{\varepsilon} d)  \times (Y,D)$$ 
such that $\displaystyle dom \ \rho_{\varepsilon}^{x} = \bar{B}(x,\varepsilon)$, 
$\displaystyle im \ \rho_{\varepsilon}^{x} = Y$, $\displaystyle (x,y) \in \rho_{\varepsilon}^{x}$ 
for any $\displaystyle \varepsilon \in (0,1]$ and 
\item[-] a "zoom modulus" $\displaystyle F: (0,1) \rightarrow [0,+\infty)$ such
that  $\displaystyle \lim_{\varepsilon \rightarrow 0} F(\varepsilon)  =  0$, 
\end{enumerate}
such that for all $\varepsilon \in (0,1)$  we have 
$\displaystyle acc(\rho_{\varepsilon}^{x}) \leq F(\varepsilon)$. 
\label{deftangsp}
\end{definition}

Using the notation proposed previously, we can write 
$\displaystyle F(\varepsilon) =  \mathcal{O}(\varepsilon)$, if there is no need 
to precisely specify a zoom modulus function.

Let us write again the definition of resolution, accuracy, precision, in the presence of
scale. The accuracy of $\displaystyle \rho_{\varepsilon}^{x}$ is defined by: 
\begin{equation}
acc(\rho_{\varepsilon}^{x}) \, = \, \sup \left\{ \mid D(y_{1}, y_{2}) - \frac{1}{\varepsilon} d(x_{1},x_{2}) \mid \mbox{
: } (x_{1},y_{1})  , \, (x_{2},y_{2}) \in \rho_{\varepsilon}^{x} \right\}
\label{acc2}
\end{equation}
 The resolution of $\displaystyle \rho_{\varepsilon}^{x}$ at  $z \in Y$ is 
\begin{equation}
res(\rho_{\varepsilon}^{x})(z) \, = \, \frac{1}{\varepsilon} \, \sup \left\{ d(x_{1},x_{2}) \mbox{
: } (x_{1},z) \in \rho_{\varepsilon}^{x} \, , \, (x_{2},z) \in \rho_{\varepsilon}^{x} \right\}
\label{resy2}
\end{equation}
and the resolution of $\displaystyle \rho_{\varepsilon}^{x}$ is given by: 
\begin{equation}
res(\rho_{\varepsilon}^{x}) \, = \,  \sup \left\{ res(\rho_{\varepsilon}^{x})(y) \mbox{
: } y \in Y \right\}
\label{res2}
\end{equation}
 The precision of $\rho_{\varepsilon}^{x}$ at $\displaystyle u \in
\bar{B}(x,\varepsilon)$ is 
\begin{equation}
prec(\rho_{\varepsilon}^{x})(u) \, = \, \sup \left\{ D(y_{1},y_{2}) \mbox{
: } (u,y_{1}) \in \rho_{\varepsilon}^{x} \, , \, (u,y_{2}) \in \rho_{\varepsilon}^{x} \right\}
\label{precx2}
\end{equation}
and the precision of $\displaystyle \rho_{\varepsilon}^{x}$ is given by: 
\begin{equation}
prec(\rho_{\varepsilon}^{x}) \, = \, \sup \left\{ prec(\rho_{\varepsilon}^{x})(u) \mbox{
: } u \in \bar{B}(x,\varepsilon) \right\}
\label{prec2}
\end{equation}

If $(Y,D, y)$  represents the (pointed unit ball in the) metric 
tangent space at $x \in X$ of $(X,d)$ and $\displaystyle \rho_{\varepsilon}^{x}$ is the
sequence of maps at smaller and smaller scale, then we have: 

\begin{equation}
 \sup \left\{ \mid D(y_{1}, y_{2}) - \frac{1}{\varepsilon} d(x_{1},x_{2}) \mid \mbox{
: } (x_{1},y_{1}) , \, (x_{2},y_{2}) \in \rho_{\varepsilon}^{x}
\right\} \, = \, \mathcal{O}(\varepsilon)
\label{acclim}
\end{equation}

\begin{equation}
 \sup \left\{  D(y_{1}, y_{2})  \mbox{
: } (u,y_{1}) \in \rho_{\varepsilon}^{x} \, , \, (u,y_{2}) \in \rho_{\varepsilon}^{x} \, ,
\,  u \in \bar{B}(x,\varepsilon)
\right\} \, = \, \mathcal{O}(\varepsilon)
\label{preclim}
\end{equation}

\begin{equation}
  \sup \left\{  d(x_{1},x_{2}) \mbox{
: } (x_{1},z) \in \rho_{\varepsilon}^{x} \, , \, (x_{2},z) \in \rho_{\varepsilon}^{x} \, , \, z
\in Y 
\right\} \, = \, \varepsilon \, \mathcal{O}(\varepsilon)
\label{reslim}
\end{equation}

Of course, relation (\ref{acclim}) implies the other two, but it is interesting to 
notice the mechanism of rescaling.

\subsection{Scale stability. Viewpoint stability}

\vspace{.5cm}

\centerline{\includegraphics[angle=270, width=0.7\textwidth]{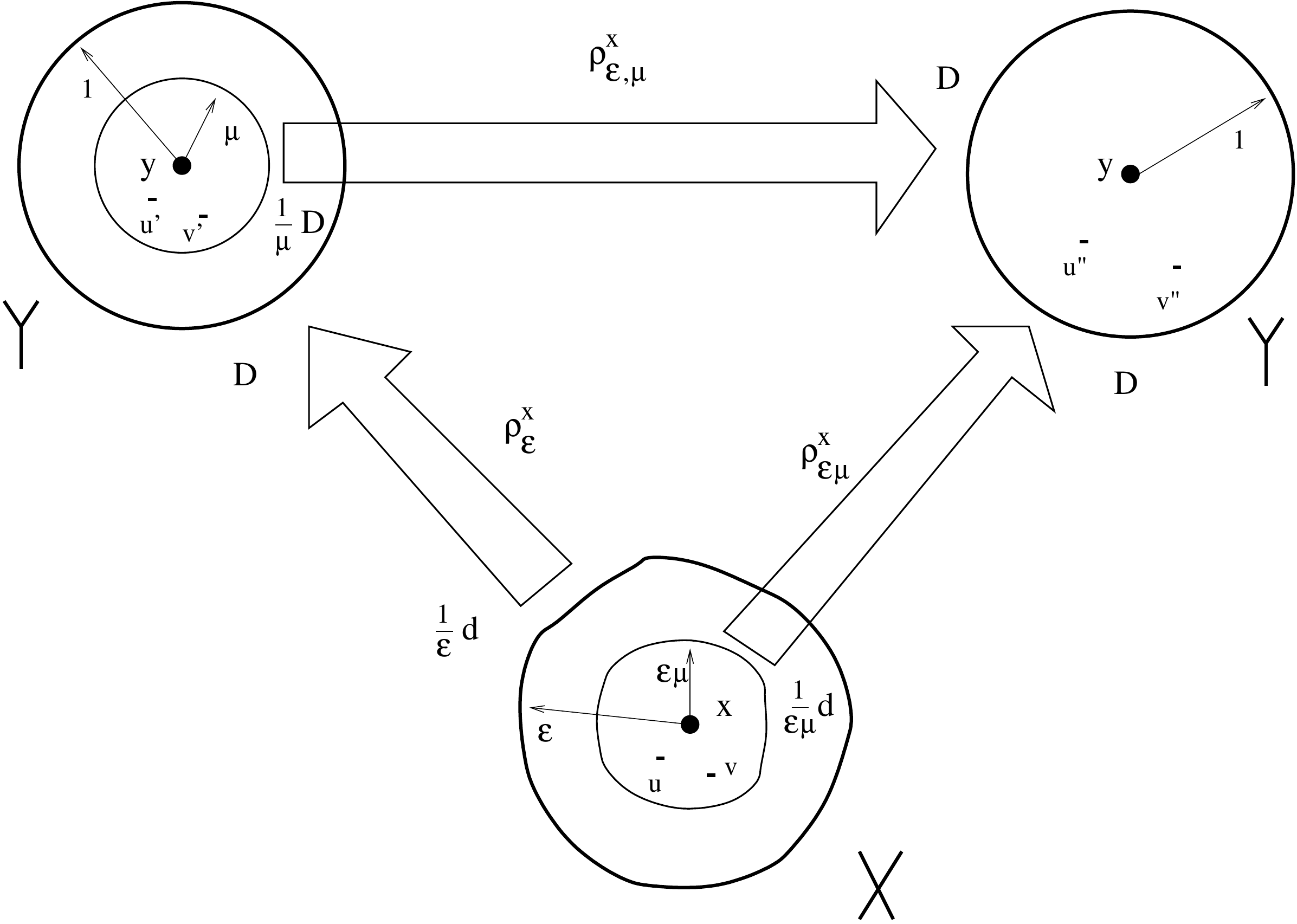}}

\vspace{.5cm}

I shall suppose further that there is a metric tangent space at $x \in X$ and I
shall work with a zoom sequence of maps described in definition \ref{deftangsp}. 

Let $\varepsilon, \mu \in (0,1)$ be two scales. Suppose we have the maps 
of the territory $X$, around $x \in X$, 
at scales $\varepsilon$ and $\varepsilon \mu$, 
$$\displaystyle 
\rho_{\varepsilon}^{x} \subset \bar{B}(x,\varepsilon) \times \bar{B}(y,1)$$ 
$$\displaystyle \rho_{\varepsilon \mu}^{x} \subset \bar{B}(x,\varepsilon 
\mu) \times \bar{B}(y,1)$$
made into the tangent space at $x$, $\displaystyle (\bar{B}(y,1), D)$. 
The ball $\displaystyle \bar{B}(x, \varepsilon \mu) \subset X$ has then 
two maps. These maps are at different scales: the first is done at scale $\varepsilon$, 
the second is done at scale $\varepsilon \mu$. 

What are the differences between these two maps? We could find out by defining
a new map 
\begin{equation}
\rho^{x}_{\varepsilon,\mu} \, = \, \left\{ (u',u") \in 
\bar{B}(y,\mu) \times \bar{B}(y,1) \mbox{ : } \right. 
\label{relchart}
\end{equation}
$$\left. \exists u \in \bar{B}(x,\varepsilon \mu) \, (u,u') \in \rho_{\varepsilon}^{x} \,
, \, (u, u") \in \rho^{x}_{\varepsilon \mu} \right\}$$ 
and measuring its accuracy, with respect to the distances $\displaystyle
\frac{1}{\mu} D$ (on the domain) and $\displaystyle D$ (on the image).

Let us consider $\displaystyle (u,u'), (v,v') \in \rho_{\varepsilon}^{x}$ and  
$\displaystyle (u,u"), (v,v") \in \rho_{\varepsilon \mu}^{x}$ such that 
$\displaystyle (u', u"), (v',v")  \in \rho^{x}_{\varepsilon,\mu}$. Then: 
$$\mid D(u",v") - \frac{1}{\mu} D(u',v') \mid \, \leq \, \mid 
\frac{1}{\mu} D(u',v') - \frac{1}{\varepsilon \mu} d(u,v) \mid \, + \, 
\mid \frac{1}{\varepsilon \mu} d(u,v) -  D(u",v") \mid$$ 
We have therefore an estimate for the accuracy of the map 
$\displaystyle \rho^{x}_{\varepsilon,\mu}$, coming from estimate (\ref{acclim})
applied for $\displaystyle \rho_{\varepsilon}^{x}$ and 
$\displaystyle \rho_{\varepsilon \mu}^{x}$: 
\begin{equation}
acc(\rho^{x}_{\varepsilon,\mu}) \, \leq \, \frac{1}{\mu} 
\mathcal{O}(\varepsilon) + \mathcal{O}(\varepsilon \mu)
\label{eqcascade}
\end{equation}
This explains the cascading of errors phenomenon, namely, for fixed $\mu$, as 
$\varepsilon$ goes to $0$ the accuracy of the map $\displaystyle
\rho^{x}_{\varepsilon,\mu}$ becomes smaller and smaller, meaning that the 
maps of the ball $\displaystyle \bar{B}(x, \varepsilon \mu) \subset X$ at the 
scales $\varepsilon, \varepsilon \mu$ (properly rescaled) are  more and more
alike. On the contrary, for fixed $\varepsilon$, as $\mu$ goes to $0$, the bound
on the accuracy becomes bigger and bigger, meaning that by using only the map 
at scale $\varepsilon$,  magnifications of  a smaller
scale region of this map may be less accurate than the map of this smaller 
region done at the smaller scale.

I shall add a supplementary hypothesis to the one 
concerning the existence of the metric tangent space. It is somehow natural to
suppose that as $\varepsilon$ converges to $0$ the map $\displaystyle
\rho^{x}_{\varepsilon,\mu}$ converges to a map $\displaystyle
\bar{\rho}^{x}_{\mu}$. This is described further. 

\vspace{.5cm}

\centerline{\includegraphics[angle=270, width=0.7\textwidth]{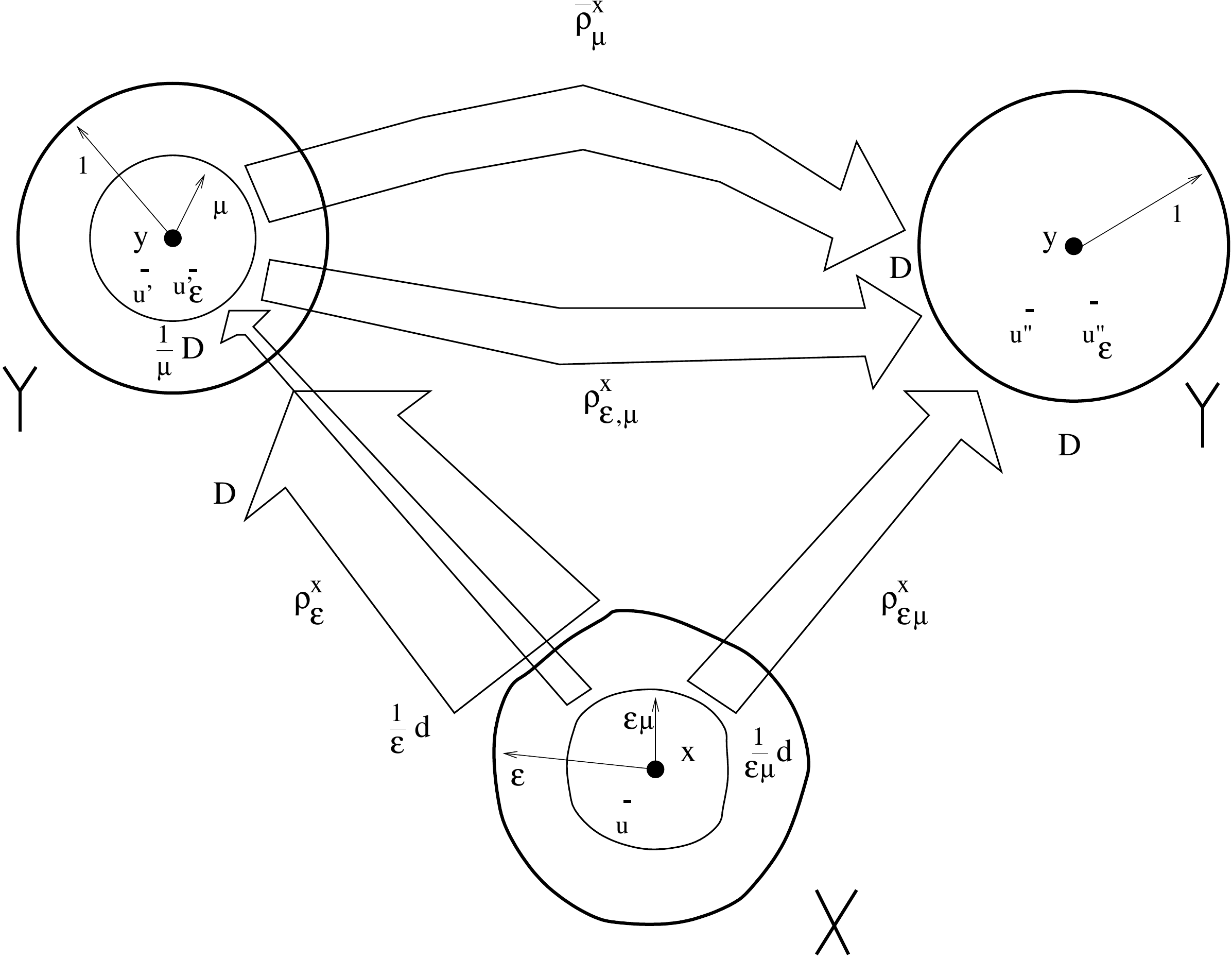}}

\vspace{.5cm}

\begin{definition}
Let the zoom sequence $\displaystyle \rho^{x}_{\varepsilon}$ be as in definition
\ref{deftangsp} and for given $\mu \in (0,1)$, the map $\displaystyle
\rho^{x}_{\varepsilon, \mu}$ be defined as in (\ref{relchart}). We say that 
the zoom sequence $\displaystyle \rho^{x}_{\varepsilon}$ is  scale stable at 
scale $\mu$ if there is a relation $\displaystyle \bar{\rho}^{x}_{\mu} \subset 
\bar{B}(y,\mu) \times \bar{B}(y,1)$ such that the Haussorff distance between 
$\displaystyle \rho^{x}_{\varepsilon, \mu}$ and $\displaystyle
\bar{\rho}^{x}_{\mu}$, in the metric space $\displaystyle \bar{B}(y,\mu) \times
\bar{B}(y,1)$ 
with the distance 
$$D_{\mu} \left( (u',u"), (v',v") \right) \, = \, \frac{1}{\mu} D(u',v') +
D(u",v") $$
can be estimated as: 
$$D_{\mu}^{Hausdorff}\left( \rho^{x}_{\varepsilon, \mu} , \bar{\rho}^{x}_{\mu}
\right)  \,
\leq \, F_{\mu}(\varepsilon)$$ 
with $\displaystyle F_{\mu}(\varepsilon) = \mathcal{O}_{\mu}(\varepsilon)$.
Such a function $\displaystyle F_{\mu}(\cdot)$ is called a scale stability modulus of 
the zoom sequence  $\displaystyle \rho^{x}_{\varepsilon}$.
\label{defstable}
\end{definition}

This means that for any $\displaystyle (u',u") \in \bar{\rho}^{x}_{\mu}$ there
is a sequence $\displaystyle (u'_{\varepsilon}, u"_{\varepsilon}) \in
\rho^{x}_{\varepsilon, \mu}$ such that 
$$\lim_{\varepsilon \rightarrow 0} u'_{\varepsilon} \, = \, u' \quad
\lim_{\varepsilon \rightarrow 0} u"_{\varepsilon} \, = \, u"$$

\begin{proposition}
If there is a scale stable zoom sequence   $\displaystyle \rho^{x}_{\varepsilon}$ 
 as in definitions \ref{deftangsp} and \ref{defstable} then the space 
 $(Y,D)$ is self-similar in a neighbourhood of point $y \in Y$, namely for any 
 $\displaystyle (u',u"), (v',v") \in \bar{\rho}^{x}_{\mu}$ we have: 
 $$D(u",v") \, = \, \frac{1}{\mu} D(u',v')$$
 In particular $\displaystyle \bar{\rho}^{x}_{\mu}$ is the graph of a function 
 (the precision and resolution are respectively equal to $0$). 
 \label{pstable}
 \end{proposition}
 
 \paragraph{Proof.} 
 Indeed, for any $\varepsilon \in (0,1)$ let us consider $\displaystyle 
 (u'_{\varepsilon}, u"_{\varepsilon}), (v'_{\varepsilon},v"_{\varepsilon}) 
 \in \rho^{x}_{\varepsilon,\mu}$ such that 
 $$\frac{1}{\mu} D(u',u'_{\varepsilon}) + D(u",u"_{\varepsilon}) \, \leq \,
 \mathcal{O}_{\mu}(\varepsilon)$$  
 $$\frac{1}{\mu} D(v',v'_{\varepsilon}) + D(v",v"_{\varepsilon}) \, \leq \,
 \mathcal{O}_{\mu}(\varepsilon)$$  
 Then we get the following inequality, using also the cascading of errors
 inequality (\ref{eqcascade}),  
 $$\mid D(u",v") - \frac{1}{\mu} D(u',v') \mid \, \leq \, 2
 \mathcal{O}_{\mu}(\varepsilon) + \frac{1}{\mu} 
\mathcal{O}(\varepsilon) + \mathcal{O}(\varepsilon \mu)$$ 
We pass with $\varepsilon$ to $0$ in order to obtain the conclusion. \hfill
$\square$

\vspace{.5cm}

Instead of changing the scale (i.e. understanding the scale stability of 
the zoom sequence), we could explore what happens when we change the point of
view. 

\vspace{.5cm}

\centerline{\includegraphics[angle=270, width=0.7\textwidth]{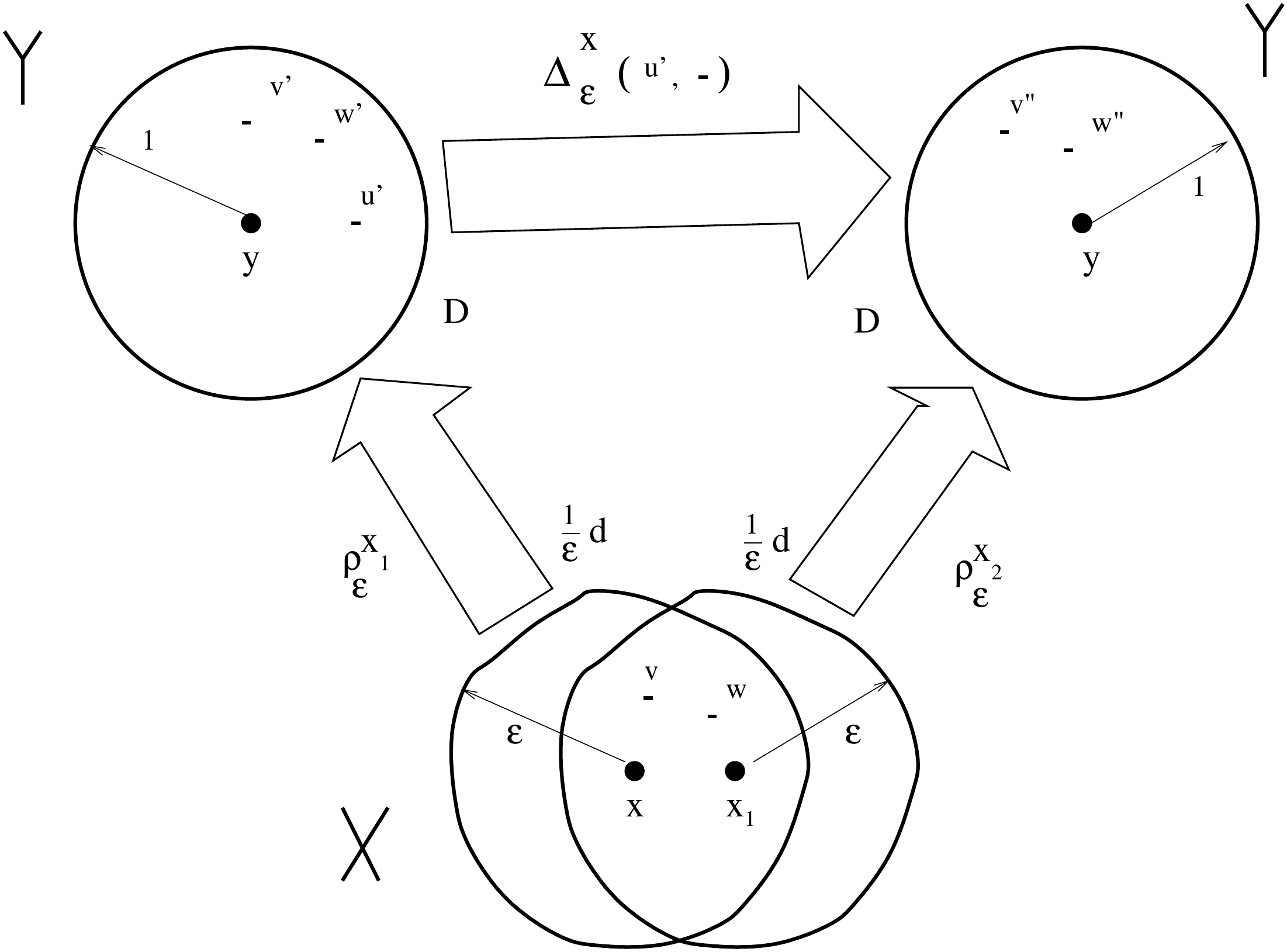}}

\vspace{.5cm}

This time we have a zoom sequence, a scale $\varepsilon \in (0,1)$ and two 
points: $\displaystyle x \in X$ and $\displaystyle u' \in \bar{B}(y,1)$. To the 
point $u'$ from the map space $Y$ corresponds a point $\displaystyle 
x_{1} \in \bar{B}(x,\varepsilon)$ such that 
$$(x_{1}, u') \in \rho_{\varepsilon}^{x}$$

The points $\displaystyle x, x_{1}$ are neighbours, in the sense that 
$\displaystyle d(x,x_{1}) < \varepsilon$. The points of $X$ which are in the 
intersection 
$$\bar{B}(x,\varepsilon) \cap \bar{B}(x_{1},\varepsilon)$$ 
are represented by both maps, $\displaystyle \rho^{x}_{\varepsilon}$ and 
$\displaystyle \rho^{x_{1}}_{\varepsilon}$. These maps are different; the
relative map between them is defined as: 
\begin{equation}
\Delta^{x}_{\varepsilon}(u', \cdot) \, = \, 
\left\{ (v',v") \in \bar{B}(y,1) \mbox{ : } \exists v \in
\bar{B}(x,\varepsilon) \cap \bar{B}(x_{1},\varepsilon) \right.
\label{diffrel}
\end{equation}
$$ \left. \quad \quad \, 
(v,v') \in \rho^{x}_{\varepsilon} \, , \, (v,v") \in 
\rho^{x_{1}}_{\varepsilon} \right\}$$
and it is called "difference at scale $\varepsilon$, from $\displaystyle x$  
to $\displaystyle x_{1}$, as seen from $\displaystyle u'$". 

The viewpoint stability of the zoom sequence is expressed as the scale
stability: the zoom sequence is stable if the difference at scale 
$\varepsilon$ converges in the sense of Hausdorff distance, as $\varepsilon$
goes to $0$. 

\begin{definition}
Let the zoom sequence $\displaystyle \rho^{x}_{\varepsilon}$ be as in definition
\ref{deftangsp} and for any $\displaystyle u' \in \bar{B}(y,1)$, 
the map $\displaystyle \Delta^{x}_{\varepsilon}(u', \cdot)$ be defined as 
in (\ref{diffrel}). The zoom sequence $\displaystyle \rho^{x}_{\varepsilon}$ 
is  viewpoint stable if there is a relation $\displaystyle \Delta^{x}(u', \cdot) \subset 
\bar{B}(y,1) \times \bar{B}(y,1)$ such that the Haussorff distance
can be estimated as: 
$$D_{\mu}^{Hausdorff}\left( \Delta^{x}_{\varepsilon}(u', \cdot) , 
\Delta^{x}(u', \cdot)
\right)  \,
\leq \, F_{diff}(\varepsilon)$$ 
with $\displaystyle F_{diff}(\varepsilon) = \mathcal{O}(\varepsilon)$.
Such a function $\displaystyle F_{diff}(\cdot)$ is called a viewpoint stability modulus of 
the zoom sequence  $\displaystyle \rho^{x}_{\varepsilon}$.
\label{defstablev}
\end{definition}

There is a proposition analoguous with proposition \ref{pstable}, stating that 
the difference relation $\displaystyle \Delta^{x}(u', \cdot)$ is the graph of 
an isometry of $(Y,D)$. 

\subsection{Foveal maps} 
The following proposition shows that if we have a scale stable zoom sequence of maps  $\displaystyle \rho^{x}_{\varepsilon}$ 
 as in definitions \ref{deftangsp} and \ref{defstable} then we can improve every
 member of the sequence such that all maps from the new zoom sequence have better
 accuracy near the "center" of the map $x \in X$, which justifies the name 
 "foveal maps".
 
\begin{definition}
Let   $\displaystyle \rho^{x}_{\varepsilon}$ be a scale stable zoom sequence. We define for any $\varepsilon \in (0,1)$ the 
 $\mu$-foveal    map $\displaystyle \phi^{x}_{\varepsilon}$  made of all pairs 
$\displaystyle (u,u') \in  \bar{B}(x,\varepsilon) \times \bar{B}(y,1)$ such 
that
\begin{enumerate}
\item[-] if $\displaystyle u \in \bar{B}(x,\varepsilon \mu) $ then
$\displaystyle (u, \bar{\rho}_{\mu}^{x}(u')) \in \rho_{\varepsilon \mu}^{x}$, 
\item[-] or else $\displaystyle (u,u') \in \rho^{x}_{\varepsilon}$.
\end{enumerate}
\label{deffoveal}
\end{definition}

\begin{proposition}
Let   $\displaystyle \rho^{x}_{\varepsilon}$ be a scale stable zoom sequence with
associated zoom modulus $F(\cdot)$ and scale stability modulus $\displaystyle
F_{\mu}(\cdot)$. The sequence of $\mu$-foveal maps $\displaystyle
\phi^{x}_{\varepsilon}$  is then a scale stable 
zoom sequence with zoom modulus $\displaystyle F(\cdot) + \mu F_{\mu}(\cdot)$. 
Moreover, the accuracy of the restricted foveal map $\displaystyle
\phi^{x}_{\varepsilon} \cap \left(\bar{B}(x,\varepsilon \mu) \times 
\bar{B}(y,\mu)\right)$ 
is bounded by $\mu F(\varepsilon \mu)$, therefore the right hand side term in
the cascading of errors inequality (\ref{eqcascade}), applied for the restricted
foveal map, can be improved to $2 F(\varepsilon \mu)$.
\label{pfoveal}
\end{proposition}

\paragraph{Proof.}
Let $\displaystyle u \in \bar{B}(x,\varepsilon \mu)$. Then there are 
$\displaystyle u', u'_{\varepsilon} \in \bar{B}(y,\mu)$ and 
$\displaystyle u",u"_{\varepsilon} \in \bar{B}(y,1)$ such that 
$\displaystyle (u,u')\in \phi^{x}_{\varepsilon}$, $\displaystyle 
(u,u") \in \rho^{x}_{\varepsilon} \mu)$, $\displaystyle (u',u") 
\in \bar{\rho}^{x}_{\mu}$, $\displaystyle (u'_{\varepsilon},u"_{\varepsilon})
\in \rho^{x}_{\varepsilon,\mu}$ and 
$$ \frac{1}{\mu} D(u',u'_{\varepsilon}) + D(u",u"_{\varepsilon}) \, \leq \,
F_{\mu}(\varepsilon)$$

Let $\displaystyle u,v \in \bar{B}(x,\varepsilon \mu)$ and $u',v' \in 
\bar{B}(y,\mu)$ such that $\displaystyle (u,u'), (v,v') \in
\phi^{x}_{\varepsilon}$. According to the definition of 
$\displaystyle \phi^{x}_{\varepsilon}$, it follows that there are 
uniquely defined $u",v" \in \bar{B}(y,1)$ such that $\displaystyle 
(u,u"), (v,v") \in \rho^{x}_{\varepsilon \mu}$ and $\displaystyle 
(u', u"), (v',v") \in \bar{\rho}^{x}_{\mu}$. We then have: 
$$\mid \frac{1}{\varepsilon} d(u,v) - D(u',v') \mid \, = \, $$
$$= \, \mid \frac{1}{\varepsilon} d(u,v) - \mu D(u",v") \mid \, = \,$$
$$= \, \mu \mid \frac{1}{\varepsilon \mu} d(u,v) - D(u",v") \mid \, \leq \, 
\mu F(\varepsilon \mu)$$
Thus we proved that the accuracy of the restricted foveal map 
$$\displaystyle
\phi^{x}_{\varepsilon} \cap \left(\bar{B}(x,\varepsilon \mu) \times
\bar{B}(y,\mu)\right)$$
is bounded by $\mu F(\varepsilon \mu)$: 
\begin{equation}
\mid \frac{1}{\varepsilon} d(u,v) - D(u',v') \mid \, \leq \, 
\mu F(\varepsilon \mu)
\label{nd1}
\end{equation}
 If $\displaystyle u, v \in \bar{B}(x,\varepsilon) \setminus \bar{B}(x,\mu)$
  and  $\displaystyle (u,u'), (v,v') \in \phi^{x}_{\varepsilon}$ then  $\displaystyle (u,u'), (v,v') \in
\rho^{x}_{\varepsilon}$, therefore 
$$\mid \frac{1}{\varepsilon} d(u,v) - D(u',v') \mid \, \leq \, F(\varepsilon)$$
Suppose now that $\displaystyle (u,u'), (v,v') \in
\phi^{x}_{\varepsilon}$ and $\displaystyle u \in \bar{B}(x,\varepsilon \mu)$ but
$\displaystyle v 
\in \bar{B}(x,\varepsilon) \setminus \bar{B}(x,\mu)$. 
 We have then: 
$$\mid \frac{1}{\varepsilon} d(u,v) - D(u',v') \mid \, \leq \, $$
$$\leq \, \mid \frac{1}{\varepsilon} d(u,v) - D(u'_{\varepsilon},v') \mid + 
D(u',u'_{\varepsilon}) \, \leq \, F(\varepsilon) + \mu F_{\mu}(\varepsilon)$$ 
We proved that the sequence of $\mu$-foveal maps $\displaystyle
\phi^{x}_{\varepsilon}$  is a  
zoom sequence with zoom modulus $\displaystyle F(\cdot) + \mu F_{\mu}(\cdot)$. 

In order to prove that the sequence is scale stable, we have to compute 
$\displaystyle \phi^{x}_{\varepsilon, \mu}$, graphically shown in the next
figure.

\vspace{.5cm}

\centerline{\includegraphics[angle=270, width=0.7\textwidth]{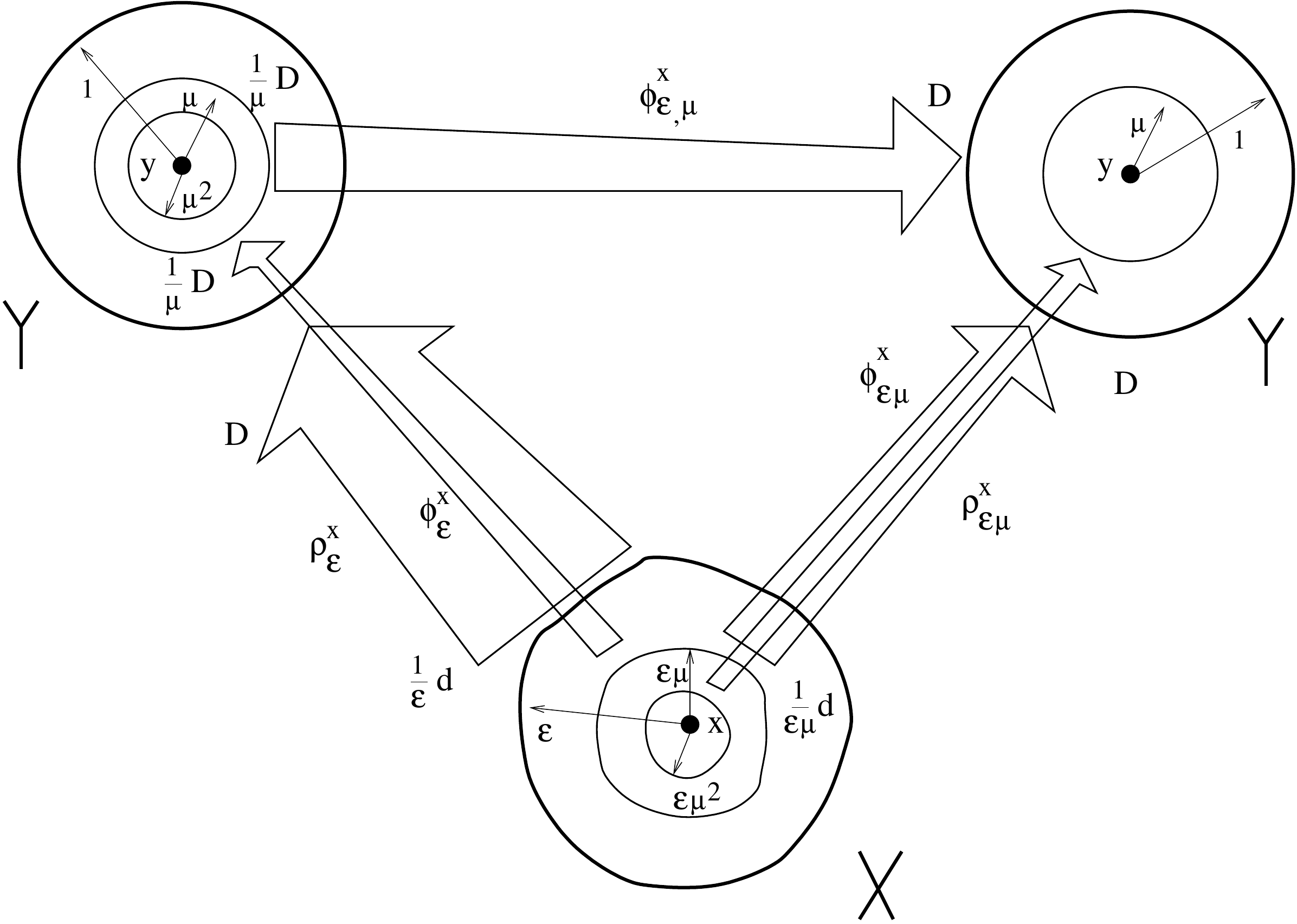}}

\vspace{.5cm}

We see that $\displaystyle (u',u") \in \phi^{x}_{\varepsilon, \mu}$ implies 
that $\displaystyle (u',u") \in \rho^{x}_{\varepsilon, \mu}$ or 
$\displaystyle (u',u") \in \rho^{x}_{\varepsilon \mu, \mu}$. From here we deduce
that the sequence of foveal maps is scale stable and that 
$$\varepsilon \mapsto \max \left\{ F_{\mu}(\varepsilon) , \mu
F_{\mu}(\varepsilon \mu) \right\}$$ 
is a scale stability modulus for the foveal sequence.

The improvement of the right hand side for the cascading of errors
inequality (\ref{eqcascade}), applied for the restricted
foveal map is then straightforward if we use (\ref{nd1}).  \hfill $\square$

\vspace{.5cm}

\section{Dilation structures}

From  definition \ref{deffoveal} we see that 
\begin{equation}
\bar{\rho}^{x}_{\mu} \circ \phi^{x}_{\varepsilon} \, = \, 
\rho^{x}_{\varepsilon \mu}
\label{pre1param}
\end{equation}

Remark that if the $\mu$-foveal map $\displaystyle \phi^{x}_{\varepsilon}$ 
coincides with the chart $\displaystyle \rho^{x}_{\varepsilon}$ for every 
$\varepsilon$ (that is, if 
the zoom sequence $\displaystyle \rho^{x}_{\varepsilon}$ is already so good 
that it cannot be improved by the construction of foveal maps), then relation 
(\ref{pre1param}) becomes 

\begin{equation}
\bar{\rho}^{x}_{\mu} \circ \phi^{x}_{\varepsilon} \, = \, 
\phi^{x}_{\varepsilon \mu}
\label{pre2param}
\end{equation}

By proposition \ref{pstable}, it follows that $\mu$-foveal map at scale 
$\varepsilon \mu$ is just a $1/\mu$ dilation of a part of the $\mu$-foveal map at scale 
$\varepsilon$. 

An idealization of these "perfect", stable zoom sequences which cannot be
improved by the $\mu$-foveal map construction for any $\mu \in (0,1)$, are
dilation structures. 

There are several further assumptions, which clearly amount to yet other 
idealizations. These are the following:
\begin{enumerate}
\item[-] the "map is the territory assumption", namely $\displaystyle 
Y = U(x)$, the "map space" 
is included in $X$, the "territory", and $y = x$. 
\item[-] "functions instead relations", that is the perfect stable zoom
sequences $\displaystyle \rho^{x}_{\varepsilon} = \phi^{x}_{\varepsilon}$ 
are graphs of functions, called dilations. That means: 
$$\rho^{x}_{\varepsilon} \subset \, \left\{ (\delta^{x}_{\varepsilon} u' , u') 
\mbox{ : } u' \in Y = V_{\varepsilon}(x) \right\}$$
\item[-] "hidden uniformity", that is: in order to pass to the limit in various 
situations, we could choose the  zoom modulus  and stability
modulus to not depend on $x \in X$. This innocuous assumption is the least
obvious, but necessary one. 
\end{enumerate}

With these idealizations in force, remember that we want our dilations to form a stable zoom sequence and 
we want also the subtler viewpoint stability, which consists in 
being  able to change the point of view in a coherent way, 
as the scale goes to zero. These are the axioms of a dilation structure. 

We shall use here a slightly particular version of dilation structures. 
For the general definition of a dilation structure see \cite{buligadil1}. 
More about this, as well as about length dilation structures, see  
\cite{buligadil3}.

\begin{definition}
Let $(X,d)$ be a complete metric space such that for any $x  \in X$ the 
closed ball $\bar{B}(x,3)$ is compact. A dilation structure $(X,d, \delta)$ 
over $(X,d)$ is the assignment to any $x \in X$  and $\varepsilon \in (0,+\infty)$ 
of a  homeomorphism, defined as: if 
$\displaystyle   \varepsilon \in (0, 1]$ then  $\displaystyle 
 \delta^{x}_{\varepsilon} : U(x)
\rightarrow V_{\varepsilon}(x)$, else 
$\displaystyle  \delta^{x}_{\varepsilon} : 
W_{\varepsilon}(x) \rightarrow U(x)$,  with the following properties.  
\begin{enumerate}
\item[{\bf A0.}]  For any $x \in X$ the sets $ \displaystyle U(x), V_{\varepsilon}(x), 
W_{\varepsilon}(x)$ are open neighbourhoods of $x$.  There are   $1<A<B$ such that for any $x \in X$  and any 
$\varepsilon \in (0,1)$ we have: 
$$\displaystyle  B_{d}(x, \varepsilon) \subset \delta^{x}_{\varepsilon}  B_{d}(x, A) 
\subset V_{\varepsilon}(x) \subset $$ 
$$\subset W_{\varepsilon^{-1}}(x) \subset \delta_{\varepsilon}^{x}  B_{d}(x, B)$$
Moreover for any compact set $K \subset X$ there are $R=R(K) > 0$ and 
$\displaystyle \varepsilon_{0}= \varepsilon(K) \in (0,1)$  such that  
for all $\displaystyle u,v \in \bar{B}_{d}(x,R)$ and all 
$\displaystyle \varepsilon  \in (0,\varepsilon_{0})$,  we have 
$\displaystyle 
\delta_{\varepsilon}^{x} v \in W_{\varepsilon^{-1}}(
\delta^{x}_{\varepsilon}u)$. 

\item[{\bf A1.}]  For any $x \in X$ 
$\displaystyle  \delta^{x}_{\varepsilon} x = x $ and $\displaystyle \delta^{x}_{1} = id$. 
Consider  the closure $\displaystyle Cl(dom \, \delta)$ of the set 
$$ dom \, \delta = \left\{ (\varepsilon, x, y) \in (0,+\infty) \times X 
\times X \mbox{ : } \right.$$ 
$$\left. \mbox{ if } \varepsilon \leq 1 \mbox{ then } y 
\in U(x) \,
\, , \mbox{  else } y \in W_{\varepsilon}(x) \right\} $$ 
seen in  $[0,+\infty) \times X \times X$ endowed with
 the product topology. The function $\displaystyle \delta : dom \, \delta 
\rightarrow  X$,  $\displaystyle \delta (\varepsilon,  x, y)  = 
\delta^{x}_{\varepsilon} y$ is continuous, admits a continuous extension 
over $\displaystyle Cl(dom \, \delta)$ and we have 
$\displaystyle \lim_{\varepsilon\rightarrow 0} \delta_{\varepsilon}^{x} y \, =
\, x$. 

\item[{\bf A2.}] For any  $x, \in X$, $\displaystyle \varepsilon, \mu \in (0,+\infty)$
 and $\displaystyle u \in U(x)$, whenever one of the sides are well defined
   we have the equality 
$\displaystyle  \delta_{\varepsilon}^{x} \delta_{\mu}^{x} u  =
\delta_{\varepsilon \mu}^{x} u$.

\item[{\bf A3.}]  For any $x$ there is a distance  function $\displaystyle (u,v) \mapsto d^{x}(u,v)$, defined for any $u,v$ in the closed ball (in distance d) $\displaystyle 
\bar{B}(x,A)$, such that uniformly with respect to $x$ in compact set we have
the limit: 
$$\lim_{\varepsilon \rightarrow 0} \quad \sup  \left\{  \mid 
\frac{1}{\varepsilon} d(\delta^{x}_{\varepsilon} u, \delta^{x}_{\varepsilon} v) \ - \ d^{x}(u,v) \mid \mbox{ :  } u,v \in \bar{B}_{d}(x,A)\right\} \ =  \ 0$$

\item[{\bf A4.}] Let us define 
$\displaystyle \Delta^{x}_{\varepsilon}(u,v) =
\delta_{\varepsilon^{-1}}^{\delta^{x}_{\varepsilon} u} \delta^{x}_{\varepsilon} v$. 
Then we have the limit, uniformly with respect to $x, u, v$ in compact set,  
$$\lim_{\varepsilon \rightarrow 0}  \Delta^{x}_{\varepsilon}(u,v) =  \Delta^{x}(u, v)  $$
 
\end{enumerate}
\label{defweakstrong}
\end{definition}

 It is algebraically straightforward to transport a dilation structure: given 
$(X,d,  \delta)$ a dilation structure and $f: X \rightarrow Z$ a uniformly
continuous homeomorphism from $X$ (as a topological space) to another topological
space $Z$ (actually more than a topological space, it should be a space endowed
with an uniformity), we can define the transport of $(X,d,\delta)$ by $f$ as 
the dilation structure $(Z,f*d, f*\delta)$. The distance $f*d$ is defined as 
$$\left(f*d\right) (u,v) = d(f(u),f(v))$$
which is a true distance, because we supposed $f$ to be a homeomorphism. 
For any $u,v \in X$ and $\varepsilon > 0$, we define the new dilation based 
at $f(u) \in Z$, of coefficient $\varepsilon$, applied to $f(v) \in Z$ as
$$\left(f*\delta\right)^{f(u)}_{\varepsilon} f(v) \, = \, f \left(
\delta^{u}_{\varepsilon} v \right)$$
It is easy to check that this is indeed a dilation structure. 

In particular we may consider to transport a dilation structure by one of its 
dilations. Visually, this corresponds to transporting the atlas representing 
a dilation structure on $X$ to a neighbourhood of one of its points. It is like 
a scale reduction of the whole territory $(X,d)$ to a smaller set. 

Inversely, we may transport the (restriction of the) dilation structure 
$(X,d,\delta)$ from $\displaystyle V_{\varepsilon}(x)$ to $U(x)$, by using 
$\displaystyle \delta^{x}_{\varepsilon^{-1}}$ as the transport function $f$. 
This is like a magnification of the "infinitesimal neighbourhood" 
$\displaystyle V_{\varepsilon}(x)$. (This neighbourhood is infinitesimal in the
sense that we may consider $\varepsilon$ as a variable, going to $0$ when
needed. Thus, instead of one neighbourhood $\displaystyle V_{\varepsilon}(x)$, 
there is a sequence of them, smaller and smaller). 

This is useful, because it allows us to make "infinitesimal statements", i.e. 
statements concerning this sequence of magnifications, as $\varepsilon
\rightarrow 0$. 

Let us compute then the magnified dilation structure. We should also rescale 
the distance on $\displaystyle V_{\varepsilon}(x)$ by a factor $1/\varepsilon$. 
Let us compute this magnified dilation structure: 
\begin{enumerate}
\item[-] the space is $U(x)$
\item[-] for any $u,v \in U(x)$ the (transported) distance between them is 
$$d^{x}_{\varepsilon}(u,v) \, = \, \frac{1}{\varepsilon}
d(\delta^{x}_{\varepsilon} u , \delta^{x}_{\varepsilon} v ) $$
\item[-] for any $u,v \in U(x)$ and scale parameter $\mu \in (0,1)$ 
(we could take $\mu > 0$ but then we have to be careful with the domains and
codomains of these new dilations), the transported dilation based at $u$, of 
coefficient $\mu$, applied to $v$, is 
\begin{equation}
 \delta^{x}_{\varepsilon^{-1}} \, \delta_{\varepsilon}^{\delta^{x}_{\varepsilon} u} \, 
\delta^{x}_{\varepsilon} v
\label{firstchora}
\end{equation}
\end{enumerate}

It is visible that working with such combinations of dilations becomes quickly
difficult. This is one of the reasons of looking for more graphical notations. 

Here is the definition of "linearity" and "selfsimilarity" for dilation structures. 

\begin{definition}
Let $X$, $Y$ be metric spaces endowed with dilation structures. 
A function $f: X \rightarrow Y$ is linear if and only if it is a morphism 
of  dilation structures: for any $u,v \in X$ and any $\varepsilon \in
\Gamma$  
$$f(\delta^{u}_{\varepsilon} v) \, = \, \delta^{f(u)}_{\varepsilon} f(v)$$ 
which is also a Lipschitz map from $X$ to $Y$ as metric spaces. 

A dilation structure $(X,d,\delta)$ is $(x,\mu)$ self-similar (for a 
$x \in X$ and $\mu\in \Gamma$, different from $1$, the neutral element
of $\Gamma = (0,+\infty)$) if the dilation $\displaystyle f = \delta^{x}_{\mu}$ is 
linear from $(X,d,\delta)$ to itself and moreover for any $u,v \in X$ we have 
$$d(\delta^{x}_{\mu} u, \delta^{x}_{\mu} v) \, = \, \mu \, d(u,v)$$

A dilation structure is linear if it is self-similar with respect to any $x \in
X$ and $\mu \in \Gamma$. 
\label{dline}
\end{definition}

\begin{definition}
Let  $(X,d, \delta)$ be a dilation structure. A property 
$$\displaystyle \mathcal{P}(x_{1},x_{2},x_{3}, ...)$$ 
is true  for $\displaystyle x_{1}, x_{2}, x_{3}, ... \in X$  sufficiently 
close if for any compact, non empty set $K \subset X$, there
is a positive constant $C(K)> 0$ such that $\displaystyle \mathcal{P}(x_{1},x_{2},
x_{3}, ...)$ is true for any $\displaystyle x_{1},x_{2},
x_{3}, ... \in K$ with $\displaystyle d(x_{i}, x_{j}) \leq C(K)$.
\end{definition}

For a dilation structure the metric tangent spaces  have the algebraic
structure of a normed group with dilations.

 We shall work further with local groups, which are spaces endowed with 
a locally defined  operation which satisfies  the conditions of a uniform group. 
 See section 3.3 \cite{buligadil1} for details about the 
definition of local groups.

\subsection{Normed conical groups}

This name has been introduced in section 8.2 \cite{buligadil1}, but these
objects appear more or less in the same form under the name "contractible
group" or "homogeneous group". Essentially these are groups endowed with a 
family of "dilations". They were also studied  in section 4 
\cite{buligadil2}. 

In the following general definition appear a topological 
commutative group $\Gamma$ endowed with a continuous morphism $\nu: \Gamma \rightarrow (0, +\infty)$ from $\Gamma$ 
to the group $(0, +\infty)$ with multiplication.  
The morphism $\nu$ induces an invariant topological filter on $\Gamma$ (other names for such an invariant filter are "absolute" or "end").  The convergence 
of a variable $\varepsilon \in \Gamma$ to this filter is denoted by $\varepsilon \rightarrow 0$ and 
it means simply $\nu(\varepsilon) \rightarrow 0$ in $\mathbb{R}$. 

Particular, interesting examples of pairs $(\Gamma, \nu)$ are: $(0, +\infty)$ with identity, which is 
the case interesting for this paper, $\displaystyle \mathbb{C}^{*}$ with the modulus of complex numbers, 
or $\mathbb{N}$ (with addition) with the exponential, which is relevant for the case of normed 
contractible groups, section 4.3 \cite{buligadil2}. 

\begin{definition}
A  normed group with dilations $(G, \delta, \| \cdot \|)$ is a local  group 
$G$  with  a local action of $\Gamma$ (denoted by $\delta$), on $G$ such that
\begin{enumerate}
\item[H0.] the limit  $\displaystyle \lim_{\varepsilon \rightarrow 0}
\delta_{\varepsilon} x  =  e$ exists and is uniform with respect to $x$ in a compact neighbourhood of the identity $e$.
\item[H1.] the limit
$\displaystyle \beta(x,y)  =  \lim_{\varepsilon \rightarrow 0} \delta_{\varepsilon}^{-1}
\left((\delta_{\varepsilon}x) (\delta_{\varepsilon}y ) \right)$
is well defined in a compact neighbourhood of $e$ and the limit is uniform with
respect to $x$, $y$.
\item[H2.] the following relation holds: 
$\displaystyle \lim_{\varepsilon \rightarrow 0} \delta_{\varepsilon}^{-1}
\left( ( \delta_{\varepsilon}x)^{-1}\right)  =  x^{-1}$, 
where the limit from the left hand side exists in a neighbourhood $U \subset G$  of 
$e$ and is uniform with respect to $x \in U$.
\end{enumerate}

Moreover the group is endowed with a continuous norm
function $\displaystyle \|\cdot \| : G \rightarrow \mathbb{R}$ which satisfies
(locally, in a neighbourhood  of the neutral element $e$) the properties:
 \begin{enumerate}
 \item[(a)] for any $x$ we have $\| x\| \geq 0$; if $\| x\| = 0$ then $x=e$,
 \item[(b)] for any $x,y$ we have $\|xy\| \leq \|x\| + \|y\|$,
 \item[(c)] for any $x$ we have $\displaystyle \| x^{-1}\| = \|x\|$,
 \item[(d)] the limit
$\displaystyle \lim_{\varepsilon \rightarrow 0} \frac{1}{\nu(\varepsilon)} \| \delta_{\varepsilon} x \| = \| x\|^{N}$
 exists, is uniform with respect to $x$ in compact set,
 \item[(e)] if $\displaystyle \| x\|^{N} = 0$ then $x=e$.
  \end{enumerate}
  \label{dnco}
  \end{definition}

\begin{theorem}
(Thm. 15 \cite{buligadil1}) Let $(G, \delta, \| \cdot \|)$ be  a locally compact  normed local group with dilations. 
Then $(G, d, \delta)$ is a  dilation structure, where the dilations 
$\delta$ and the distance $d$ are defined by: 
$\displaystyle  \delta^{x}_{\varepsilon} u = x \delta_{\varepsilon} ( x^{-1}u)  \quad , \quad
d(x,y) = \| x^{-1}y\|$. 

Moreover $(G, d, \delta)$ is linear, in the sense of definition \ref{dline}. 
\label{tgrd}
\end{theorem}

\begin{definition}
A normed conical group $N$ is a normed  group with dilations  such that
for any $\varepsilon \in \Gamma$  the dilation
 $\delta_{\varepsilon}$ is a group morphism  and such that for any $\varepsilon >0$
  $\displaystyle  \| \delta_{\varepsilon} x \| = \nu(\varepsilon) \| x \|$.
\end{definition}

A normed conical group is the infinitesimal version of a normed group with
dilations (\cite{buligadil1} proposition 2).

\begin{proposition}
Let $(G, \delta, \| \cdot \|)$ be  a locally compact  normed local group with
dilations. Then  $\displaystyle (G,\beta, \delta, \| \cdot \|^{N})$ is a locally compact, 
local normed conical group, with operation $\beta$,  dilations $\delta$ and homogeneous norm $\displaystyle \| \cdot \|^{N}$.
\label{here3.4}
\end{proposition}

\subsection{Tangent bundle of a dilation structure}
\label{induced}

The most important metric and algebraic first order  properties of a dilation 
structure are  presented here as condensed statements, available in full length as theorems 7, 8,
10 in \cite{buligadil1}.

\begin{theorem}
Let $(X,d,\delta)$ be a  dilation structure. Then the metric space $(X,d)$ 
admits a metric tangent space at $x$, for any point $x\in X$. 
More precisely we have  the following limit: 
$$\lim_{\varepsilon \rightarrow 0} \ \frac{1}{\varepsilon} \sup \left\{  \mid d(u,v) - d^{x}(u,v) \mid \mbox{ : } d(x,u) \leq \varepsilon \ , \ d(x,v) \leq \varepsilon \right\} \ = \ 0 \ .$$
\label{thcone}
\end{theorem}

\begin{theorem}
  If $(X,d,\delta)$ is a  dilation structure 
 then for any $x \in X$ the triple $\displaystyle (U^{x}, \delta^{x}, d^{x})$ is
 a locally compact normed conical group, with operation 
 $\displaystyle \Sigma^{x}(\cdot, \cdot)$, neutral element $x$ and inverse 
 $\displaystyle inv^{x}(y) = \Delta^{x}(y, x)$.  
\label{tgene}
\end{theorem}

The conical group $\displaystyle (U(x), \Sigma^{x}, \delta^{x})$ can be seen as the tangent space 
of $(X,d, \delta)$ at $x$. We shall  denote it by  
$\displaystyle T_{x} (X, d, \delta) =  (U(x), \Sigma^{x}, \delta^{x})$, or by $\displaystyle T_{x} X$ if 
$(d,\delta)$ are clear from the context.

The following proposition is  corollary 
6.3 from \cite{buligadil2}, which gives a more precise
description of the conical group 
$\displaystyle (U(x), \Sigma^{x}, \delta^{x})$.  In the proof of that corollary
there is a gap pointed by S. Vodopyanov, namely that  
 Siebert' proposition 5.4 \cite{siebert}, which 
is true for conical groups (in our language), is used for local  
conical groups. Fortunately, this gap was filled  by the 
theorem 1.1 \cite{recent}, which states that a locally compact, 
locally connected, contractible (with  Siebert' wording) group is 
locally isomorphic to a contractive Lie group.

\begin{proposition}
Let $(X,d,\delta)$ be a  dilation structure. 
Then for any $x \in X$ the local group
 $\displaystyle (U(x), \Sigma^{x})$ is locally a simply connected Lie group
 whose Lie algebra admits a positive graduation (a homogeneous group), given by
 the eigenspaces of $\displaystyle \delta^{x}_{\varepsilon}$ for an 
 arbitrary $\varepsilon \in (0,1)$.
\label{cor63}
\end{proposition}

There is a bijection between linear (in the sense of definition 
\ref{dline})   dilation structures and normed conical groups. Any 
normed conical group induces a linear  dilation structure, by theorem 
\ref{tgrd}. Conversely, we have the following result (see  theorem 6.1 
\cite{buligairq} for a more general statement).

\begin{theorem}
Let $(G,d, \delta)$ be a linear  dilation structure. Then, with the 
notations from  theorem \ref{tgene}, for any 
$x \in G$, the dilation structure $(U(x), d, \delta)$ coincides with 
the dilation structure of the conical group $\displaystyle (U(x), \Sigma^{x},
\delta^{x})$. 
\label{pgroudlin}
\end{theorem}

\subsection{Differentiability with respect to dilation structures}

For any  dilation structure or there is an associated  notion  of 
differentiability (section 7.2 \cite{buligadil1}). 
For defining differentiability with respect to dilation structures we need 
first the  definition of a morphism of conical groups. 

\begin{definition}
 Let $(N,\delta)$ and $(M,\bar{\delta})$ be two  conical groups. A function $f:N\rightarrow M$ is a conical group morphism if $f$ is a group morphism and for any $\varepsilon>0$ and $u\in N$ we have 
 $\displaystyle f(\delta_{\varepsilon} u) = \bar{\delta}_{\varepsilon} f(u)$. 
\label{defmorph}
\end{definition}

The definition of the derivative, or differential,  with respect to dilations 
structures is a straightforward generalization of  the 
definition of the Pansu derivative  \cite{pansu}.

 \begin{definition}
 Let $(X, d, \delta)$ and $(Y, \overline{d}, \overline{\delta})$ be two 
  dilation structures   and $f:X \rightarrow Y$ be a continuous function. The function $f$ is differentiable in $x$ if there exists a 
 conical group morphism  $\displaystyle D \, f(x):T_{x}X\rightarrow T_{f(x)}Y$, defined on a neighbourhood of $x$ with values in  a neighbourhood  of $f(x)$ such that 
\begin{equation}
\lim_{\varepsilon \rightarrow 0} \sup \left\{  \frac{1}{\varepsilon} \overline{d} \left( f\left( \delta^{x}_{\varepsilon} u\right) ,  \overline{\delta}^{f(x)}_{\varepsilon} D \, f(x)  (u) \right) \mbox{ : } d(x,u) \leq \varepsilon \right\}Ê  = 0 , 
\label{edefdif}
\end{equation}
The morphism $\displaystyle D \, f(x) $ is called the derivative, or differential,  of $f$ at $x$.
\label{defdiffer}
\end{definition}

The definition also makes sense if the function $f$ is defined on a 
open subset of $(X,d)$.

\end{document}